\pdfoutput=1
\RequirePackage{ifpdf}
\ifpdf 
\documentclass[pdftex]{sigma}
\else
\documentclass{sigma}
\fi

\numberwithin{equation}{section}

\newtheorem{Theorem}{Theorem}[section]

\newtheorem{Proposition}[Theorem]{Proposition}
{\theoremstyle{definition}
\newtheorem{Definition}[Theorem]{Definition}

}

\usepackage[all,cmtip]{xy}

\newcommand*\pfqskip{8mu}
\catcode`,\active
\newcommand*\pfq{\begingroup
        \catcode`\,\active
        \def ,{\mskip\pfqskip\relax}%
        \dopfq
}
\catcode`\,12
\def\dopfq#1#2#3#4#5{%
        {}_{#1}\phi_{#2}\bigg(\genfrac..{0pt}{}{#3}{#4}\Big|#5\biggr)%
        \endgroup
}
\newcommand*\pFqskip{8mu}
\catcode`,\active
\newcommand*\pFq{\begingroup
        \catcode`\,\active
        \def ,{\mskip\pFqskip\relax}%
        \dopFq
}
\catcode`\,12
\def\dopFq#1#2#3#4#5{%
        {}_{#1}F_{#2}\biggl[\genfrac..{0pt}{}{#3}{#4};#5\biggr]%
        \endgroup
}

\begin{document}

\allowdisplaybreaks

\renewcommand{\PaperNumber}{038}

\FirstPageHeading

\ShortArticleName{A ``Continuous'' Limit of the Complementary Bannai--Ito Polynomials}

\ArticleName{A ``Continuous'' Limit of the Complementary\\
Bannai--Ito Polynomials: Chihara Polynomials}

\Author{Vincent X.~GENEST~$^\dag$, Luc VINET~$^\dag$ and Alexei ZHEDANOV~$^\ddag$}

\AuthorNameForHeading{V.X.~Genest, L.~Vinet and A.~Zhedanov}

\Address{$^\dag$~Centre de Recherches Math\'ematiques, Universit\'e de Montr\'eal,\\
\hphantom{$^\dag$}~C.P.~6128, Succ.
Centre-Ville, Montr\'eal, QC, Canada, H3C 3J7}

\EmailD{\href{mailto:genestvi@crm.umontreal.ca}{genestvi@crm.umontreal.ca},
\href{mailto:vinetl@crm.umontreal.ca}{vinetl@crm.umontreal.ca}}

\Address{$^\ddag$~Donetsk Institute for Physics and Technology, Donetsk 83114, Ukraine}
\EmailD{\href{mailto:zhedanov@yahoo.com}{zhedanov@yahoo.com}}

\ArticleDates{Received December 23, 2013, in f\/inal form March 24, 2014; Published online March 30, 2014}

\Abstract{A novel family of $-1$ orthogonal polynomials called the Chihara polynomials is characterized.
The polynomials are obtained from a~``continuous'' limit of the complementary Bannai--Ito polynomials, which are the
kernel partners of the Bannai--Ito polynomials.
The three-term recurrence relation and the explicit expression in terms of Gauss hypergeometric functions are obtained
through a~limit process.
A one-parameter family of second-order dif\/ferential Dunkl operators having these polynomials as eigenfunctions is also
exhibited.
The quadratic algebra with involution encoding this bispectrality is obtained.
The orthogonality measure is derived in two dif\/ferent ways: by using Chihara's method for kernel polynomials and, by
obtaining the symmetry factor for the one-parameter family of Dunkl operators.
It is shown that the polynomials are related to the big $-1$ Jacobi polynomials by a~Christof\/fel transformation and that
they can be obtained from the big $q$-Jacobi by a~$q\rightarrow -1$ limit.
The generalized Gegenbauer/Hermite polynomials are respectively seen to be special/limiting cases of the Chihara
polynomials.
A one-parameter extension of the generalized Hermite polynomials is proposed.}

\Keywords{Bannai--Ito polynomials; Dunkl operators; orthogonal polynomials; quadratic algebras}

\Classification{33C45}

\section{Introduction}

One of the recent advances in the theory of orthogonal polynomials is the characterization of~$-1$ orthogonal
polynomials~\cite{Genest-2013-02-1,
Genest-2012,Tsujimoto-2012-03,Tsujimoto-2013-03-01,Vinet-2011-02,Vinet-2011-01,Vinet-2012-05}.
The distinguishing property of these polynomials is that they are eigenfunctions of Dunkl-type operators which involve ref\/lections.
They also correspond to $q\rightarrow -1$ limits of certain $q$-polynomials of the Askey tableau.
The~$-1$ polynomials should be organized in a~tableau complementing the latter.
Sitting atop this $-1$ tableau would be the Bannai--Ito polynomials (BI) and their kernel partners, the complementary
Bannai--Ito polynomials (CBI).
Both families depend on four real parameters, satisfy a~discrete/f\/inite orthogonality relation and correspond to
a~(dif\/ferent) $q\rightarrow -1$ limit of the Askey--Wilson polynomials~\cite{Askey-1985}.
The BI polynomials are eigenfunctions of a~\emph{first-order} Dunkl dif\/ference operator whereas the CBI polynomials are
eigenfunctions of a~\emph{second-order} Dunkl dif\/ference operator.
It should be noted that the polynomials of the $-1$ scheme do not all have the same type of bispectral properties in
distinction with what is observed when $q\rightarrow 1$ because the second-order $q$-dif\/ference equations of the basic
polynomials of the Askey scheme do not always exist in certain $q\rightarrow -1$ limits.
In this paper, a~novel family of $-1$ orthogonal polynomials stemming from a~``continuous'' limit of the complementary
Bannai--Ito polynomials will be studied and characterized.
Its members will be called Chihara polynomials.

The Bannai--Ito polynomials, written $B_{n}(x;\rho_1,\rho_2,r_1,r_2)$ in the notation of~\cite{Genest-2013-02-1}, were
f\/irst identif\/ied by Bannai and Ito themselves in their classif\/ication~\cite{Bannai-1984} of orthogonal polynomials
satisfying the Leonard duality property~\cite{Leonard-1982-07}; they were also seen to correspond to a~$q\rightarrow -1$
limit of the $q$-Racah polynomials~\cite{Bannai-1984}.
A signif\/icant step in the characterization of the BI polynomials was made in~\cite{Tsujimoto-2012-03} where it was
recognized that the polynomials $B_{n}(x)$ are eigenfunctions of the most general (self-adjoint) f\/irst-order Dunkl shift
operator which stabilizes polynomials of a~given degree, i.e.
\begin{gather}
\label{BI-OP}
\mathcal{L}=\left[\frac{(x-\rho_1)(x-\rho_2)}{2x}\right](\mathbb{I}-R)+\left[\frac{(x-r_1+1/2)(x-r_2+1/2)}{2x+1}\right](T^{+}R-\mathbb{I}),
\end{gather}
where $T^{+}f(x)=f(x+1)$ is the shift operator, $Rf(x)=f(-x)$ is the ref\/lection operator and where $\mathbb{I}$ stands
for the identity.
In the same paper~\cite{Tsujimoto-2012-03}, it was also shown that the BI polynomials correspond to a~$q\rightarrow -1$
limit of the Askey--Wilson polynomials and that the operator~\eqref{BI-OP} can be obtained from the Askey--Wilson
operator in this limit.
An important limiting case of the BI polynomials is found by considering the ``continuous'' limit, which is obtained
upon writing
\begin{gather}
\label{Continuous-Limit}
x\rightarrow \frac{x}{h},
\qquad
\rho_1=\frac{a_1}{h}+b_1,
\qquad
\rho_2=\frac{a_2}{h}+b_2,
\qquad
r_1=\frac{a_1}{h},
\qquad
r_2=\frac{a_2}{h},
\end{gather}
and taking $h\rightarrow 0$.
In this limit, the operator~\eqref{BI-OP} becomes, after a~rescaling of the variable $x$, the most general
(self-adjoint) f\/irst-order \emph{differential} Dunkl operator which preserves the space of polynomials of a~given
degree, i.e.
\begin{gather}
\label{-1-Jacobi-OP}
\mathcal{M}=\left[\frac{(a+b+1)x^2+(ac-b)x+c}{x^2}\right](R-\mathbb{I})+\left[\frac{2(x-1)(x+c)}{x}\right]\partial_{x}R.
\end{gather}
The polynomial eigenfunctions of~\eqref{-1-Jacobi-OP} have been identif\/ied in~\cite{Vinet-2011-02,Vinet-2012-05} as the
big $-1$ Jacobi polynomials $J_{n}(x;a,b,c)$ introduced in~\cite{Vinet-2012-05}.
Alternatively, one can obtain the polynomials $J_{n}(x;a,b,c)$ by directly applying the limit~\eqref{Continuous-Limit}
to the BI polynomials.
The big $-1$ Jacobi polynomials satisfy a~continuous orthogonality relation on the interval $[-1,-c]\cup[1,c]$.
They also correspond to a~$q\rightarrow -1$ limit of the big $q$-Jacobi polynomials, an observation which was f\/irst used
to derive their properties in~\cite{Vinet-2012-05}.
It is known moreover (see for example~\cite{Koekoek-2010}) that the big $q$-Jacobi polynomials can be obtained from the
Askey--Wilson polynomials using a~limiting procedure similar to~\eqref{Continuous-Limit}.
Hence the relationships between the Askey--Wilson, big $q$-Jacobi, Bannai--Ito and big $-1$ Jacobi polynomials can be
expressed diagrammatically as follows
\begin{gather}
\label{First-Diagram}
\begin{gathered}
\xymatrix
@C=.15\linewidth
@R=.12\linewidth
{
*+{
\begin{aligned}
&\text{Askey--Wilson}\\&p_{n}(x;a,b,c,d\,|\,q)
\end{aligned}}
\ar@{->}[r]^-{q\rightarrow -1}
\ar@{->}[d]^-{a\rightarrow 0}_-{x\rightarrow x/2a}
&
*+{\begin{aligned}
&\text{Bannai--Ito}\\
&B_{n}(x;\rho_1,\rho_2,r_1,r_2)
\end{aligned}}
\ar@{->}[d]^-{h\rightarrow 0}_-{x\rightarrow x/h}
\\
*+{
\begin{aligned}
&\text{big } q\text{-Jacobi} \\
&P_{n}(x;\alpha,\beta,\gamma\,|\, q)
\end{aligned}
}
\ar@{->}[r]^-{q\rightarrow -1}
&
*+{
\begin{aligned}
&\text{big } {-}1\text{ Jacobi}\\
& J_{n}(x;a,b,c)
\end{aligned}
}
}
\end{gathered}
\end{gather}
where the notation of~\cite{Koekoek-2010} was used for the Askey--Wilson $p_{n}(x;a,b,c,d\,|\,q)$ and the big
$q$-Jacobi polynomials $P_{n}(x;\alpha,\beta,\gamma\,|\,q)$.

In this paper, we shall be concerned with the continuous limit~\eqref{Continuous-Limit} of the complementary Bannai--Ito
polynomials $I_{n}(x;\rho_1,\rho_2,r_1,r_2)$~\cite{Genest-2013-02-1,Tsujimoto-2012-03}.
The polynomials $C_{n}(x;\alpha,\beta,\gamma)$ arising in this limit shall be referred to as Chihara polynomials since
they have been introduced by T.~Chihara in~\cite{Chihara-1964} (up to a~parameter redef\/inition).
They depend on three real parameters.
Using the limit, the recurrence relation and the explicit expression for the polynomials $C_{n}(x;\alpha,\beta,\gamma)$
in terms of Gauss hypergeometric functions will be obtained from that of the CBI polynomials.
The second-order dif\/ferential Dunkl operator having the Chihara polynomials as eigenfunctions will also be given.
The corresponding bispectrality property will be used to construct the algebraic structure behind the Chihara
polynomials: a~quadratic Jacobi algebra~\cite{Granovskii-1992-07} supplemented with an involution.
The weight function for the Chihara polynomials will be constructed in two dif\/ferent ways: on the one hand using
Chihara's method for kernel polynomials~\cite{Chihara-1964} and on the other hand by solving a~Pearson-type
equation~\cite{Koekoek-2010}.
This measure will be def\/ined on the union of two disjoint intervals.
The Chihara polynomials $C_{n}(x;\alpha,\beta,\gamma)$ will also be seen to correspond to a~$q\rightarrow -1$ limit of
the big $q$-Jacobi polynomials that is dif\/ferent from the one leading to the big $-1$ Jacobi.
In analogy with~\eqref{First-Diagram}, the following relationships shall be established:
\begin{gather*}
\begin{gathered}
\xymatrix
@C=.15\linewidth
@R=.08\linewidth
{
*+{
\begin{aligned}
&\text{Askey--Wilson}\\&p_{n}(x;a,b,c,d\,|\,q)
\end{aligned}}
\ar@{->}[r]^-{q\rightarrow -1}
\ar@{->}[d]^-{a\rightarrow 0}_-{x\rightarrow x/2a}
&
*+{\begin{aligned}
&\text{complementary BI}\\
&I_{n}(x;\rho_1,\rho_2,r_1,r_2)
\end{aligned}}
\ar@{->}[d]^-{h\rightarrow 0}_-{x\rightarrow x/h}
\\
*+{
\begin{aligned}
&\text{big } q\text{-Jacobi} \\
&P_{n}(x;\alpha,\beta,\gamma\,|\, q)
\end{aligned}
}
\ar@{->}[r]^-{q\rightarrow -1}
&
*+{
\begin{aligned}
&\text{Chihara}\\
& C_{n}(x;\alpha,\beta,\gamma)
\end{aligned}
}
}
\end{gathered}
\end{gather*}
Since the CBI polynomials are obtained from the BI polynomials by the Christof\/fel transform~\cite{Chihara-2011} (and
vice-versa using the Geronimus transform~\cite{Geronimus-1947}), it will be shown that the following relations relating
the Chihara to the big $-1$ Jacobi polynomials hold:
\begin{gather*}
\begin{gathered}
\xymatrix
@C=.15\linewidth
@R=.1\linewidth
{
*+{
\begin{aligned}
&\text{Bannai--Ito}\\&B_{n}(x;\rho_1,\rho_2,r_1,r_2)
\end{aligned}}
\ar@<.5ex>[r]^-{\text{Christoffel}}
\ar@{->}[d]^-{h\rightarrow 0}_-{x\rightarrow x/h}
&
*+{\begin{aligned}
&\text{complementary BI}\\
&I_{n}(x;\rho_1,\rho_2,r_1,r_2)
\end{aligned}}
\ar@<.5ex>[l]^-{\text{Geronimus}}
\ar@{->}[d]^-{h\rightarrow 0}_-{x\rightarrow x/h}
\\
*+{
\begin{aligned}
&\text{big } {-}1\text{ Jacobi} \\
&J_{n}(x;a,b,c)
\end{aligned}
}
\ar@<.5ex>[r]^-{\text{Christoffel}}
&
*+{
\begin{aligned}
&\text{Chihara}\\
& C_{n}(x;\alpha,\beta,\gamma)
\end{aligned}
}
\ar@<.5ex>[l]^-{\text{Geronimus}}
}
\end{gathered}
\end{gather*}
Finally, it will be observed that for $\gamma=0$, the Chihara polynomials $C_{n}(x;\alpha,\beta,\gamma)$ reduce to the
generalized Gegenbauer polynomials and that upon taking the limit $\beta\rightarrow \infty$ with $\gamma=0$, the
polynomials $C_{n}(x;\alpha,\beta,\gamma)$ go to the generalized Hermite polynomials~\cite{Chihara-2011}.
A one-parameter extension of the generalized Hermite polynomials will also be presented.

The remainder of the paper is organized straightforwardly.
In Section~\ref{Section2}, the main features of the CBI polynomials are reviewed.
In Section~\ref{Section3}, the ``continuous'' limit is used to def\/ine the Chihara polynomials and establish their basic properties.
In Section~\ref{Section4}, the operator having the Chihara polynomials as eigenfunctions is obtained and the algebraic structure
behind their bispectrality is exhibited.
In Section~\ref{Section5}, the weight function is derived and the orthogonality relation is given.
In Section~\ref{Section6}, the polynomials are related to the big $-1$ Jacobi and big $q$-Jacobi polynomials.
In Section~\ref{Section7}, limits and special cases are examined.

\section{Complementary Bannai--Ito polynomials}\label{Section2}

In this section, the main properties of the complementary Bannai--Ito polynomials, which have been obtained
in~\cite{Genest-2013-02-1, Tsujimoto-2012-03}, are reviewed.
Let $\rho_1$, $\rho_2$, $r_1$, $r_2$ be real parameters, the monic CBI polynomial $I_{n}(x;\rho_1,\rho_2,r_1,r_2)$,
denoted $I_{n}(x)$ for notational ease, are def\/ined by
\begin{gather}
I_{2n}(x)=\eta_{2n}\pFq{4}{3}{-n,n+g+1,\rho_2+x,\rho_2-x}{\rho_1+\rho_2+1,\rho_2-r_1+1/2,\rho_2-r_2+1/2}{1},
\nonumber
\\
I_{2n+1}(x)=\eta_{2n+1}(x-\rho_2)\pFq{4}{3}{-n,n+g+2,\rho_2+x+1,\rho_2-x+1}{\rho_1+\rho_2+2,\rho_2-r_1+3/2,\rho_2-r_2+3/2}{1},
\label{CBI-DEF}
\end{gather}
where $g=\rho_1+\rho_2-r_1-r_2$ and where $_{p}F_{q}$ denotes the generalized hypergeometric series~\cite{Gasper-2004}.
It is directly seen from~\eqref{CBI-DEF} that $I_{n}(x)$ is a~polynomial of degree $n$ in $x$ and that it is symmetric
with respect to the exchange of the two parameters $r_1$, $r_2$.
The coef\/f\/icients $\eta_{n}$, which ensure that the polynomials are monic (i.e.\
$I_{n}(x)=x^{n}+\mathcal{O}(x^{n-1})$), are given by
\begin{gather*}
\eta_{2n}=\frac{(\rho_1+\rho_2+1)_{n}(\rho_2-r_1+1/2)_{n}(\rho_2-r_2+1/2)_{n}}{(n+g+1)_{n}},
\\
\eta_{2n+1}=\frac{(\rho_1+\rho_2+2)_{n}(\rho_2-r_1+3/2)_{n}(\rho_2-r_2+3/2)_{n}}{(n+g+2)_{n}},
\end{gather*}
where $(a)_{n}=a(a+1)\cdots(a+n-1)$, $(a)_{0}\equiv 1$, stands for the Pochhammer symbol.
The CBI polynomials satisfy the three-term recurrence relation
\begin{gather}
\label{CBI-RECU}
xI_{n}(x)=I_{n+1}(x)+(-1)^{n}\rho_2I_{n}(x)+\tau_{n}I_{n-1}(x),
\end{gather}
subject to the initial conditions $I_{-1}(x)=0$, $I_{0}(x)=1$ and with the recurrence coef\/f\/icients
\begin{gather}
\tau_{2n}=-\frac{n(n+\rho_1-r_1+1/2)(n+\rho_1-r_2+1/2)(n-r_1-r_2)}{(2n+g)(2n+g+1)},
\nonumber
\\
\tau_{2n+1}=-\frac{(n+g+1)(n+\rho_1+\rho_2+1)(n+\rho_2-r_1+1/2)(n+\rho_2-r_2+1/2)}{(2n+g+1)(2n+g+2)}.
\label{CBI-COEF}
\end{gather}
The CBI polynomials form a~f\/inite set $\{I_{n}(x)\}_{n=0}^{N}$ of positive-def\/inite orthogonal polynomials provided that
the truncation and positivity conditions $\tau_{N+1}=0$ and $\tau_{n}>0$ hold for $n=1,\ldots, N$, where $N$ is
a~positive integer.
Under these conditions, the CBI polynomials obey the orthogonality relation
\begin{gather*}
\sum\limits_{i=0}^{N}\omega_iI_{n}(x_i)I_{m}(x_i)=h_{n}^{(N)}\delta_{nm},
\end{gather*}
where the grid points $x_i$ are of the general form
\begin{gather*}
x_i=(-1)^{i}(a+1/4+i/2)-1/4,
\qquad
\text{or}
\qquad
x_i=(-1)^{i}(b-1/4-i/2)-1/4.
\end{gather*}
The expressions for the grid points $x_i$ and for the weight function $\omega_i$ depend on the truncation condition
$\tau_{N+1}=0$, which can be realized in six dif\/ferent ways (three for each possible parity of $N$).
The explicit formulas for each case shall not be needed here and can be found in~\cite{Genest-2013-02-1}.

One of the most important properties of the complementary Bannai--Ito polynomials is their bispectrality.
Recall that a~family of orthogonal polynomials $\{P_{n}(x)\}$ is bispectral if one has an eigenvalue equation of the
form
\begin{gather*}
\mathcal{A}P_{n}(x)=\lambda_{n}P_{n}(x),
\end{gather*}
where $\mathcal{A}$ is an operator acting on the argument $x$ of the polynomials.
For the CBI polynomials, there is a~one-parameter family of eigenvalue equations~\cite{Genest-2013-02-1}
\begin{gather}
\mathcal{K}^{(\alpha)}I_{n}(x)=\Lambda_{n}^{(\alpha)}I_{n}(x),
\end{gather}
with eigenvalues $\Lambda_{n}^{(\alpha)}$
\begin{gather}
\label{CBI-EIGEN}
\Lambda_{2n}^{(\alpha)}=n(n+g+1),
\qquad
\Lambda_{2n+1}^{(\alpha)}=n(n+g+2)+\omega+\alpha,
\end{gather}
where
\begin{gather}
\omega=\rho_1(1-r_1-r_2)+r_1r_2-3(r_1+r_2)/2+5/4,
\end{gather}
and where $\alpha$ is an arbitrary parameter.
The operator $\mathcal{K}^{(\alpha)}$ is the second-order Dunkl shift operator\footnote{One should take
$\alpha\rightarrow \omega+\alpha$ in the operator obtained in~\cite{Genest-2013-02-1} to f\/ind the
expression~\eqref{CBI-DUNKL}.}
\begin{gather}
\label{CBI-DUNKL}
\mathcal{K}^{(\alpha)}=A_{x}(T^{+}-\mathbb{I})+B_{x}(T^{-}-R)+C_{x}(\mathbb{I}-R)+D_{x}(T^{+}R-\mathbb{I}),
\end{gather}
where $T^{\pm}f(x)=f(x\pm 1)$, $Rf(x)=f(-x)$ and where the coef\/f\/icients read
\begin{gather}
A_{x}=\frac{(x+\rho_1+1)(x+\rho_2+1)(x-r_1+1/2)(x-r_2+1/2)}{2(x+1)(2x+1)},
\nonumber
\\
B_{x}=\frac{(x-\rho_1-1)(x-\rho_2)(x+r_1-1/2)(x+r_2-1/2)}{2x(2x-1)},
\nonumber
\\
C_{x}=\frac{(x+\rho_1+1)(x-\rho_2)(x-r_1+1/2)(x-r_2+1/2)}{2x(2x+1)}+\frac{(\alpha-x^2)(x-\rho_2)}{2x},
\nonumber
\\
D_{x}=\frac{\rho_2(x+\rho_1+1)(x-r_1+1/2)(x-r_2+1/2)}{2x(x+1)(2x+1)}.
\end{gather}
The complementary Bannai--Ito correspond to a~$q\rightarrow -1$ limit of the Askey--Wilson
polyno\-mials~\cite{Tsujimoto-2012-03}.
Consider the Askey--Wilson polyno\-mials~\cite{Askey-1985}
\begin{gather}
\label{AW-POLY}
p_{n}(z;a,b,c,d\,|\,q)=a^{-n}(ab,ac,ad;q)_{n}\pfq{4}{3}{q^{-n},abcdq^{n-1},az,az^{-1}}{ab,ac,ad}{q;q},
\end{gather}
where ${}_p\phi_{q}$ is the generalized $q$-hypergeometric function~\cite{Gasper-2004}.
Upon considering
\begin{gather*}
a=ie^{\epsilon(2\rho_1+3/2)},
\qquad
b=-ie^{\epsilon(2\rho_2+1/2)},
\qquad
c=ie^{\epsilon(-2r_2+1/2)},
\qquad
d=ie^{\epsilon(-2r_1+1/2)},
\\
q=-e^{\epsilon},
\qquad
z=ie^{-2\epsilon(x+1/4)},
\end{gather*}
and taking the limit $\epsilon\rightarrow 0$, one f\/inds that the polynomials~\eqref{AW-POLY} converge, up to
a~normalization factor, to the CBI polynomials $I_{n}(x;\rho_1,\rho_2,r_1,r_2)$.

\section{A ``continuous'' limit to Chihara polynomials}\label{Section3}

In this section, the ``continuous'' limit of the complementary Bannai--Ito polynomials will be used to def\/ine the
Chihara polynomials and obtain the three-term recurrence relation that they satisfy.
Let $\rho_1$, $\rho_2$, $r_1$, $r_2$ be parametrized as follows
\begin{gather}
\label{CLIM}
\rho_1=\frac{a_1}{h}+b_1,
\qquad
\rho_2=\frac{a_2}{h}+b_2,
\qquad
r_1=\frac{a_1}{h},
\qquad
r_2=\frac{a_2}{h},
\end{gather}
and denote by
\begin{gather}
\label{LIM-POLY}
F_{n}^{(h)}(x)=h^{n}I_{n}(x/h)
\end{gather}
the monic polynomials obtained by replacing $x\rightarrow x/h$ in the CBI polynomials.
Upon taking $h\rightarrow 0$ in the def\/inition~\eqref{CBI-DEF} of the CBI polynomials, where~\eqref{CLIM} has been used,
one f\/inds that the limit exists and that it yields
\begin{gather}
\lim_{h\rightarrow 0}F_{2n}^{(h)}
=\frac{(a_2^2-a_1^2)^{n}(b_2+1/2)_{n}}{(n+b_1+b_2+1)_{n}}\pFq{2}{1}{-n,n+b_1+b_2+1}{b_2+1/2}{\frac{a_2^2-x^2}{a_2^2-a_1^2}},
\nonumber
\\
\lim_{h\rightarrow 0}F_{2n+1}^{(h)}
=\frac{(a_2^2-a_1^2)^{n}(b_2+3/2)_{n}}{(n+b_1+b_2+2)_{n}}(x-a_2)\pFq{2}{1}{-n,n+b_1+b_2+2}{b_2+3/2}{\frac{a_2^2-x^2}{a_2^2-a_1^2}}.
\label{FIRST-LIMIT}
\end{gather}
It is directly seen that the variable $x$ in~\eqref{FIRST-LIMIT} can be rescaled and consequently, that there is only
three independent parameters.
Assuming that $a_1^2\neq a_2^2$, we can take
\begin{gather}
\label{PARA}
x\rightarrow x\sqrt{a_1^2-a_2^2},
\qquad
\alpha=b_2-1/2,
\qquad
\beta=b_1+1/2,
\qquad
\gamma=a_2/\sqrt{a_1^2-a_2^2},
\end{gather}
to rewrite the polynomials~\eqref{FIRST-LIMIT} in terms of the three parameters $\alpha$, $\beta$ and $\gamma$.
We shall moreover assume that $\gamma$ is real.
This construction motivates the following def\/inition.
\begin{Definition}
Let $\alpha$, $\beta$ and $\gamma$ be real parameters.
The Chihara polynomials $C_{n}(x;\alpha,\beta,\gamma)$, denoted $C_{n}(x)$ for simplicity, are the monic polynomials of
degree $n$ in the variable $x$ def\/ined by
\begin{gather}
C_{2n}(x)=(-1)^{n}\frac{(\alpha+1)_{n}}{(n+\alpha+\beta+1)_{n}}\pFq{2}{1}{-n,n+\alpha+\beta+1}{\alpha+1}{x^2-\gamma^2},
\nonumber
\\
C_{2n+1}(x)=(-1)^{n}\frac{(\alpha+2)_{n}}{(n+\alpha+\beta+2)_{n}}(x-\gamma)\pFq{2}{1}{-n,n+\alpha+\beta+2}{\alpha+2}{x^2-\gamma^2}.
\label{CHIHARA-DEF}
\end{gather}
\end{Definition}
The polynomials $C_{n}(x;\alpha,\beta,\gamma)$ (up to redef\/inition of the parameters) have been considered by Chihara
in~\cite{Chihara-1964} in a~completely dif\/ferent context (see Section~\ref{Section5}).
We shall henceforth refer to the polynomials $C_{n}(x;\alpha,\beta,\gamma)$ as the Chihara polynomials.
They correspond to the continuous limit~\eqref{CLIM},~\eqref{LIM-POLY} as $h\rightarrow 0$ of the CBI polynomials with
the scaling and reparametrization~\eqref{PARA}.

Using the same limit on~\eqref{CBI-RECU} and~\eqref{CBI-COEF}, the recurrence relation satisf\/ied by the Chihara
polynomials~\eqref{CHIHARA-DEF} can readily be obtained.
\begin{Proposition}[\cite{Chihara-1964}]\label{Prop-1}  
The Chihara polynomials $C_{n}(x)$ defined by~\eqref{CHIHARA-DEF} satisfy the recurrence relation
\begin{gather}
\label{CHIHARA-RECU}
xC_{n}(x)=C_{n+1}(x)+(-1)^{n}\gamma C_{n}(x)+\sigma_{n}C_{n-1}(x),
\end{gather}
where
\begin{gather}
\label{CHIHARA-COEF}
\sigma_{2n}=\frac{n(n+\beta)}{(2n+\alpha+\beta)(2n+\alpha+\beta+1)},
\qquad\!\!
\sigma_{2n+1}=\frac{(n+\alpha+1)(n+\alpha+\beta+1)}{(2n+\alpha+\beta+1)(2n+\alpha+\beta+2)}.\!\!\!
\end{gather}
\end{Proposition}
\begin{proof}
By taking the limit~\eqref{CLIM},~\eqref{LIM-POLY} and reparametrization~\eqref{PARA}
on~\eqref{CBI-RECU},~\eqref{CBI-COEF}.
\end{proof}

As is directly checked from the recurrence coef\/f\/icients~\eqref{CHIHARA-COEF}, the positivity condition $\sigma_{n}>0$
for $n\geqslant 1$ is satisf\/ied if the parameters $\alpha$ and $\beta$ are such that
\begin{gather*}
\alpha>-1,
\qquad
\beta>-1.
\end{gather*}
By Favard's theorem~\cite{Chihara-2011}, it follows that the system of polynomials
$\{C_{n}(x;\alpha,\beta,\gamma)\}_{n=0}^{\infty}$ def\/ined by~\eqref{CHIHARA-DEF} is orthogonal with respect to some
positive measure on the real line.
This measure shall be constructed in Section~\ref{Section5}.

\section{Bispectrality of the Chihara polynomials}\label{Section4}

In this section, the operator having the Chihara polynomials as eigenfunctions is derived and the algebraic structure
behind this bispectrality, a~quadratic algebra with an involution, is exhibited.

\subsection{Bispectrality}

Consider the family of eigenvalue equations~\eqref{CBI-EIGEN} satisf\/ied by the CBI polynomials.
Upon chan\-ging the variable $x\rightarrow x/h$, the action of the operator~\eqref{CBI-DUNKL} becomes
\begin{gather*}
\mathcal{K}^{(\alpha)}f(x)=A_{x/h}\left(f(x+h)-f(x)\right)+B_{x/h}(f(x-h)-f(-x))
\\
\phantom{\mathcal{K}^{(\alpha)}f(x)=}
{}+C_{x/h}(f(x)-f(-x))+D_{x/h}(f(-x-h)-f(x)).
\end{gather*}
Using the above expression and the parametrization~\eqref{CLIM}, the limit as $h\rightarrow 0$ can be taken
in~\eqref{CBI-EIGEN} to obtain the family of eigenvalue equations satisf\/ied by the Chihara polynomials.
\begin{Proposition}
\label{Prop-2}
Let $\epsilon$ be an arbitrary parameter.
The Chihara polynomials $C_{n}(x;\alpha,\beta,\gamma)$ satisfy the one-parameter family of eigenvalue equations
\begin{gather}
\label{CHIHARA-EIGEN}
\mathcal{D}^{(\epsilon)}C_{n}(x;\alpha,\beta,\gamma)=\lambda_{n}^{(\epsilon)}C_{n}(x;\alpha,\beta,\gamma),
\end{gather}
where the eigenvalues are given by
\begin{gather}
\label{CHIHARA-EIGEN-2}
\lambda_{2n}^{(\epsilon)}=n(n+\alpha+\beta+1),
\qquad
\lambda_{2n+1}^{(\epsilon)}=n(n+\alpha+\beta+2)+\epsilon,
\end{gather}
for $n=0,1,\ldots$.
The second-order differential Dunkl operator $\mathcal{D}^{(\epsilon)}$ having the Chihara polynomials as eigenfunctions
has the expression
\begin{gather}
\label{CHIHARA-DUNKL}
\mathcal{D}^{(\epsilon)}=S_{x}\partial_{x}^2+T_{x}\partial_{x}R+U_x\partial_{x}+V_{x}(\mathbb{I}-R),
\end{gather}
where the coefficients are
\begin{gather}
S_{x}=\frac{(x^2-\gamma^2)(x^2-\gamma^2-1)}{4x^2},
\qquad
T_{x}=\frac{\gamma(x-\gamma)(x^2-\gamma^2-1)}{4x^3},
\nonumber
\\
U_{x}=\frac{\gamma(x^2-\gamma^2-1)(2\gamma-x)}{4x^3}+\frac{(x^2-\gamma^2)(\alpha+\beta+3/2)}{2x}-\frac{\alpha+1/2}{2x},
\nonumber
\\
V_{x}=\frac{\gamma(x^2-\gamma^2-1)(x-3\gamma/2)}{4x^4}-\frac{(x^2-\gamma^2)(\alpha+\beta+3/2)}{4x^2}
+\frac{\alpha+1/2}{4x^2}+\epsilon\frac{x-\gamma}{2x}.
\label{CHIHARA-DUNKL-COEF}
\end{gather}
\end{Proposition}
\begin{proof}
We obtain $\mathcal{D}^{(0)}$ f\/irst.
Consider the operator $\mathcal{K}^{(-\omega)}$.
Upon taking $x\rightarrow x/h$, the action of this operator on functions of argument $x$ can be cast in the form
\begin{gather}
\mathcal{K}^{(-\omega)}f(x)=A_{x/h}\left[f(x+h)-f(x)\right]+B_{x/h}\left[f(x-h)-f(x)\right]
\nonumber
\\
\phantom{\mathcal{K}^{(-\omega)}f(x)=}
{}+\left[B_{x/h}+C_{x/h}-D_{x/h}\right]f(x)
+D_{x/h}f(-x-h)-\big[B_{x/h}+C_{x/h}\big]f(-x).\!\!\!
\label{Dompe-2}
\end{gather}
Assuming that $f(x)$ is an analytic function, the f\/irst term of~\eqref{Dompe-2} yields
\begin{gather*}
\lim_{h\rightarrow 0}\big(A_{x/h}\left[f(x+h)-f(x)\right]+B_{x/h}\left[f(x-h)-f(x)\right]\big)
=\left[\frac{(x^2-a_1^2)(x^2-a_2^2)}{4x^2}\right]f''(x)
\\
\qquad
{}+\frac{1}{4}\left[x(2b_1+2b_2+3)+\frac{-a_2x-a_2^2(1+2b_1)-2a_1^2b_2}{x}+\frac{a_1^2a_2}{x^2}-\frac{2a_1^2a_2^2}{x^3}\right]f'(x),
\end{gather*}
where~\eqref{CLIM} has been used and where $f'(x)$ stands for the derivative with respect to the argu\-ment~$x$.
With~\eqref{PARA} this gives the term $S_{x}\partial_x^2+U_x \partial_x$ in $\mathcal{D}^{(0)}$.
Similarly, using~\eqref{CLIM}, the second term of~\eqref{Dompe-2} produces
\begin{gather*}
\lim_{h\rightarrow 0}\big(B_{x/h}+C_{x/h}-D_{x/h}\big)f(x)
\\
\qquad
=\left[\frac{3a_1^2a_2^2}{8x^4}-\frac{a_1^2a_2}{4x^3}+\frac{a_2^2b_1+a_1^2b_2}{4x^2}+\frac{a_2}{4x}-\frac{2b_1+2b_2+3}{8}\right]f(x).
\end{gather*}
With the parametrization~\eqref{PARA}, this gives the term $V_x\mathbb{I}$ in $\mathcal{D}^{(0)}$.
The third term of~\eqref{Dompe-2} gives
\begin{gather*}
\lim_{h\rightarrow 0}\big(D_{x/h}f(-x-h)-\left[B_{x/h}+C_{x/h}\right]f(-x)\big)
=\frac{a_2(x^2-a_1^2)(a_2-x)}{4x^3}f'(-x)
\\
\qquad{}
-\left[\frac{3a_1^2a_2^2}{8x^4}-\frac{a_1^2a_2}{4x^3}+\frac{a_2^2b_1+a_1^2b_2}{4x^2}+\frac{a_2}{4x}-\frac{2b_1+2b_2+3}{8}\right]f(-x).
\end{gather*}
Using~\eqref{PARA}, this gives the term $-V_xR+ T_x\partial_{x}R$ in $\mathcal{D}^{(0)}$.
The arbitrary parameter $\epsilon$ can be added to the odd part of the spectrum since the Chihara polynomials satisfy
the eigenvalue equation
\begin{gather*}
\frac{(x-\gamma)}{2x}(\mathbb{I}-R)C_{n}(x)=\rho_{n}C_{n}(x)
\qquad
\text{with}
\qquad
\rho_n=
\begin{cases}
0 & \text{if $n$ is even},
\\
1 & \text{if $n$ is odd},
\end{cases}
\end{gather*}
as can be seen directly from the explicit expression~\eqref{CHIHARA-DEF}.
This concludes the proof.
\end{proof}

\subsection{Algebraic structure}

The bispectrality property of the Chihara polynomials can be encoded algebraically.
Let~$\kappa_1$,~$\kappa_2$ and $P$ be def\/ined as follows
\begin{gather*}
\kappa_1=\mathcal{D}^{(\epsilon)},
\qquad
\kappa_2=x,
\qquad
P=R+\frac{\gamma}{x}(\mathbb{I}-R),
\end{gather*}
where $\mathcal{D}^{(\epsilon)}$ is given by~\eqref{CHIHARA-DUNKL}, $R$ is the ref\/lection operator and where $\kappa_2$
corresponds to multiplication by $x$.
It is directly checked that $P$ is an involution, which means that
\begin{gather*}
P^2={\mathbb{I}}.
\end{gather*}
Upon def\/ining a~third generator
\begin{gather*}
\kappa_3=[\kappa_1,\kappa_2],
\end{gather*}
with $[a,b]=ab-ba$, a~direct computation shows that one has the commutation relations
\begin{gather}
[\kappa_3,\kappa_2]=\frac{1}{2}\kappa_2^2+\delta_2\kappa_2^2P+2\delta_3\kappa_3P-\delta_5 P-\delta_4,
\nonumber
\\
[\kappa_1,\kappa_3]=\frac{1}{2}\{\kappa_1,\kappa_2\}-\delta_2\kappa_3 P-\delta_3\kappa_1 P+\delta_1\kappa_2-\delta_1\delta_3 P,
\label{Algebra-1}
\end{gather}
where $\{a,b\}=ab+ba$ stands for the anticommutator.
The commutation relations involving the involution $P$ are given by
\begin{gather}
\label{Algebra-2}
[\kappa_1,P]=0,
\qquad
\{\kappa_2,P\}=2\delta_3,
\qquad
\{\kappa_3,P\}=0,
\end{gather}
and the structure constants $\delta_i$, $i=1,\ldots,5$ are expressed as follows
\begin{gather*}
\delta_1=\epsilon(\alpha+\beta+1-\epsilon),
\qquad
\delta_{2}=(\alpha+\beta+3/2-2\epsilon),
\qquad
\delta_3=\gamma,
\\
\delta_4=(\gamma^2+1)/2,
\qquad
\delta_5=\gamma^2(\alpha+\beta+3/2-2\epsilon)+\alpha+1/2.
\end{gather*}
The algebra def\/ined by~\eqref{Algebra-1} and~\eqref{Algebra-2} corresponds to a~Jacobi algebra~\cite{Granovskii-1992-07}
supplemented with involutions and can be seen as a~contraction of the complementary Bannai--Ito
algebra~\cite{Genest-2013-02-1}.

\section{Orthogonality of the Chihara polynomials}\label{Section5}

In this section, we derive the orthogonality relation satisf\/ied by the Chihara polynomials in two dif\/ferent ways.
First, the weight function will be constructed directly, following a~method proposed by Chihara in~\cite{Chihara-1964}.
Second, a~Pearson-type equation will be solved for the operator~\eqref{CHIHARA-DUNKL}.
It is worth noting here that the weight function cannot be obtained from the limit process~\eqref{CLIM} as $h\rightarrow
0$.
Indeed, while the complementary Bannai--Ito polynomials $F_{n}^{(h)}(x)=h^{n}I_{n}(x/h)$ approach this limit, they no
longer form a~(f\/inite) system of orthogonal polynomials.
A similar situation occurs in the standard limit from the $q$-Racah to the big $q$-Jacobi
polynomials~\cite{Koekoek-2010} and is discussed by Koornwinder in~\cite{Koorwinder-2011-04}.

\subsection{Weight function and Chihara's method}

Our f\/irst approach to the construction of the weight function is based on the method developed by Chihara
in~\cite{Chihara-1964} to construct systems of orthogonal polynomials from a~given a~set of ortho\-go\-nal polynomials and
their kernel partners (see also~\cite{Marcellan-1997-07}).
Since the present context is rather dif\/ferent, the analysis will be taken from the start.
The main observation is that the Chihara polynomials~\eqref{CHIHARA-DEF} can be expressed in terms of the Jacobi
polynomials $P_{n}^{(\alpha,\beta)}(x)$ as follows
\begin{gather}
C_{2n}(x;\alpha,\beta,\gamma)=\frac{(-1)^{n}n!}{(n+\alpha+\beta+1)_{n}}P_{n}^{(\alpha,\beta)}(y(x)),
\nonumber
\\
C_{2n+1}(x;\alpha,\beta,\gamma)=\frac{(-1)^{n}n!(x-\gamma)}{(n+\alpha+\beta+2)_{n}}P_{n}^{(\alpha+1,\beta)}(y(x)),
\label{JACO-IDEN}
\end{gather}
where
\begin{gather*}
y(x)=1-2x^2+2\gamma^2.
\end{gather*}
The Jacobi polynomials $P_{n}^{(\alpha,\beta)}(z)$ are known~\cite{Koekoek-2010} to satisfy the orthogonality relation
\begin{gather}
\label{JACOBI-ORTHO}
\int_{-1}^{1}P_{n}^{(\alpha,\beta)}(z)P_{m}^{(\alpha,\beta)}(z)\mathrm{d}\psi^{(\alpha,\beta)}(z)=\chi_{n}^{(\alpha,\beta)}\delta_{nm},
\end{gather}
with
\begin{gather}
\label{JACOBI-NORM}
\chi_{n}^{(\alpha,\beta)}=\frac{2^{\alpha+\beta+1}}{2n+\alpha+\beta+1}\frac{\Gamma(n+\alpha+1)\Gamma(n+\beta+1)}{\Gamma(n+\alpha+\beta+1)n!},
\end{gather}
where $\Gamma(z)$ is the gamma function and where
\begin{gather}
\label{JACOBI-WEIGHT}
\mathrm{d}\psi^{(\alpha,\beta)}(z)=(1-z)^{\alpha}(1+z)^{\beta}\mathrm{d}z.
\end{gather}
The relation~\eqref{JACOBI-ORTHO} is valid provided that $\alpha>-1$, $\beta>-1$.
Since the Chihara polynomials are orthogonal (by proposition~\ref{Prop-1} and Favard's theorem) and given the
relation~\eqref{JACO-IDEN} and orthogonality relation~\eqref{JACOBI-ORTHO}, we consider the integral
\begin{gather*}
\mathcal{I}_{MN}=\int_{\mathcal{F}}C_{M}(x)C_{N}(x)\mathrm{d}\phi(x),
\end{gather*}
where the interval $\mathcal{F}=\big[-\sqrt{1+\gamma^2},-|\gamma|\big]\cup
\big[|\gamma|,\sqrt{1+\gamma^2}\big]$ corresponds to the inverse mapping of the interval $[-1,1]$ for $y(x)$ and
where $\phi(x)$ is a~distribution function.
Upon taking $M=2m$ and using~\eqref{JACO-IDEN}, one directly has (up to normalization)
\begin{gather*}
\begin{split}
& \mathcal{I}_{2m,2n}=\int_{|\gamma|}^{\sqrt{1+\gamma^2}}P_{m}^{(\alpha,\beta)}\big(y(x)\big)P_{n}^{(\alpha,\beta)}\big(y(x)\big)
\big[\mathrm{d}\phi(x)-\mathrm{d}\phi(-x)\big],
\\
& \mathcal{I}_{2m,2n+1}=\int_{|\gamma|}^{\sqrt{1+\gamma^2}}P_{m}^{(\alpha,\beta)}\big(y(x)\big)
P_{n}^{(\alpha+1,\beta)}\big(y(x)\big)\big[(x-\gamma)\mathrm{d}\phi(x)+(x+\gamma)\mathrm{d}\phi(-x)\big].
\end{split}
\end{gather*}
In order that $\mathcal{I}_{2n,2m}=\mathcal{I}_{2n,2m+1}=0$ for $n\neq m$, one must have for $|\gamma
|\leqslant x \leqslant\sqrt{1+\gamma^2}$
\begin{gather}
\mathrm{d}\phi(x)-\mathrm{d}\phi(-x)=\mathrm{d}\psi^{(\alpha,\beta)}\big(y(x)\big),
\qquad
(x-\gamma)\mathrm{d}\phi(x)+(x+\gamma)\mathrm{d}\phi(-x)=0,
\label{Conditions-1}
\end{gather}
where $\psi^{(\alpha,\beta)}\big(y(x)\big)$ is the distribution appearing in~\eqref{JACOBI-WEIGHT} with $z=y(x)$.
The common solution to the equations~\eqref{Conditions-1} is seen to be given by
\begin{gather}
\label{CHIHARA-WEIGHT}
\mathrm{d}\phi(x)=\frac{(x+\gamma)}{2|x|}\mathrm{d}\psi^{(\alpha,\beta)}\big(1-2x^2+2\gamma^2\big).
\end{gather}
It is easily verif\/ied that the condition $\mathcal{I}_{2n+1,2m+1}=0$ for $n\neq m$ holds.
Indeed, upon using~\eqref{CHIHARA-WEIGHT} one f\/inds (up to normalization)
\begin{gather*}
\mathcal{I}_{2n+1,2m+1}=\int_{|\gamma
|}^{\sqrt{1+\gamma^2}}P_{n}^{(\alpha+1,\beta)}\big(y(x)\big)P_{m}^{(\alpha+1,\beta)}\big(y(x)\big)
\big[(x-\gamma)^2\mathrm{d}\phi(x)-(x+\gamma)^2\mathrm{d}\phi(-x)\big]
\\
\phantom{\mathcal{I}_{2n+1,2m+1}}
=\int_{|\gamma|}^{\sqrt{1+\gamma^2}}P_{n}^{(\alpha+1,\beta)}\big(y(x)\big)
P_{m}^{(\alpha+1,\beta)}\big(y(x)\big)\mathrm{d}\psi^{(\alpha+1,\beta)}\big(y(x)\big)=\chi_{n}^{(\alpha+1,\beta)}
\delta_{nm},
\end{gather*}
which follows from~\eqref{JACOBI-ORTHO}.
The following result has thus been established.
\begin{Proposition}[\cite{Chihara-1964}]
Let $\alpha$, $\beta$ and $\gamma$ be real parameters such that $\alpha,\beta>-1$.
The Chihara polynomials $C_{n}(x;\alpha,\beta,\gamma)$ satisfy the orthogonality relation
\begin{gather}
\label{CHIHARA-ORTHO}
\int_{\mathcal{E}}C_{n}(x)C_{m}(x)\omega(x)\mathrm{d}x=k_{n}\delta_{nm},
\end{gather}
on the interval
$\mathcal{E}=\big[-\sqrt{1+\gamma^2},-|\gamma|\big]\cup\big[|\gamma|,\sqrt{1+\gamma^2}\big]$.
The weight function has the expression
\begin{gather}
\label{CHIHARA-WEIGHT-3}
\omega(x)=\theta(x)(x+\gamma)(x^2-\gamma^2)^{\alpha}(1+\gamma^2-x^2)^{\beta},
\end{gather}
where $\theta(x)$ is the sign function.
The normalization factor $k_{n}$ is given by
\begin{gather}
k_{2n}=\frac{\Gamma(n+\alpha+1)\Gamma(n+\beta+1)}{\Gamma(n+\alpha+\beta+1)}\frac{n!}{(2n+\alpha+\beta+1)\left[(n+\alpha+\beta+1)_{n}\right]^{2}},
\nonumber
\\
k_{2n+1}=\frac{\Gamma(n+\alpha+2)\Gamma(n+\beta+1)}{\Gamma(n+\alpha+\beta+2)}\frac{n!}{(2n+\alpha+\beta+2)\left[(n+\alpha+\beta+2)_{n}\right]^{2}}.
\label{Norm}
\end{gather}
\end{Proposition}
\begin{proof}
The proof of the orthogonality relation follows from the above considerations.
The normalization factor is obtained by comparison with that of the Jacobi polynomials~\eqref{JACOBI-NORM}.
\end{proof}

\subsection{A Pearson-type equation}

The weight function for the Chihara polynomials $C_{n}(x)$ can also be derived from their bispectral
property~\eqref{CHIHARA-EIGEN} by solving a~Pearson-type equation.
A similar approach was adopted in~\cite{Vinet-2011-02} and led to the weight function for the big $-1$ Jacobi
polynomials.
In view of the recurrence relation~\eqref{CHIHARA-RECU} satisf\/ied by the Chihara polynomials $C_{n}(x)$, it follows from
Favard's theorem that there exists a~linear functional $\sigma$ such that
\begin{gather}
\label{LIN-FUN}
\langle \sigma,C_{n}(x)C_{m}(x)\rangle=h_{n}\delta_{nm},
\end{gather}
with non-zero constants $h_{n}$.
Moreover, it follows from~\eqref{CHIHARA-EIGEN} and from the completeness of the system of polynomials $\{C_{n}(x)\}$
that the operator $\mathcal{D}^{(\epsilon)}$ def\/ined by~\eqref{CHIHARA-DUNKL} and~\eqref{CHIHARA-DUNKL-COEF} is
symmetric with respect to the functional $\sigma$, which means that
\begin{gather*}
\big\langle \sigma,\{\mathcal{D}^{(\epsilon)}V(x)\} W(x)\big\rangle
=\big\langle\sigma,V(x)\{\mathcal{D}^{(\epsilon)}W(x)\}\big\rangle,
\end{gather*}
where $V(x)$ and $W(x)$ are arbitrary polynomials.
In the positive-def\/inite case $\alpha>-1$, $\beta>-1$, one has $h_{n}>0$ and there is a~realization of~\eqref{LIN-FUN}
in terms of an integral
\begin{gather*}
\langle \sigma, C_{n}(x)C_{m}(x)\rangle=\int_{a}^{b}C_{n}(x)C_{m}(x)\mathrm{d}\sigma(x),
\end{gather*}
where $\sigma(x)$ is a~distribution function and where $a$, $b$ can be inf\/inite.
Let us consider the case where $\omega(x)=\mathrm{d}\sigma(x)/\mathrm{d}x>0$ inside the interval $[a,b]$.
In this case, the following condition must hold:
\begin{gather}
\label{SYM-COND}
\big(\omega(x)\mathcal{D}^{(\epsilon)}\big)^{*}=\omega(x) \mathcal{D}^{(\epsilon)},
\end{gather}
where $\mathcal{A}^{*}$ denotes the Lagrange adjoint operator with respect to $\mathcal{A}$.
Recall that for a~generic Dunkl dif\/ferential operator
\begin{gather*}
\mathcal{A}=\sum\limits_{k=0}^{N}A_{k}(x)\partial_{x}^{k}+\sum\limits_{\ell=0}^{N}B_{k}(x)\partial_{x}^{k}R,
\end{gather*}
where $A_{k}(x)$ and $B_{k}(x)$ are real functions, the Lagrange adjoint operator reads~\cite{Vinet-2011-02}
\begin{gather*}
\mathcal{A}^{*}=\sum\limits_{k=0}^{N}(-1)^{k}\partial_{x}^{k}A_{k}(x)+\sum\limits_{\ell=0}^{N}\partial_{x}^{k}B_{k}(-x)R.
\end{gather*}
These formulas assume that the interval of orthogonality is necessarily symmetric.
Let us now derive directly the expression for the weight function $\omega(x)$ from the condition~\eqref{SYM-COND}.
Assuming $\epsilon\in \mathbb{R}$, the Lagrange adjoint of $\mathcal{D}^{(\epsilon)}$ reads
\begin{gather*}
\big[\mathcal{D}^{(\epsilon)}\big]^{*}=\partial_{x}^2S_{x}+\partial_{x}T_{-x}R-\partial_{x} U_{x}-V_{-x}R+V_{x}\mathbb{I},
\end{gather*}
where the coef\/f\/icients are given by~\eqref{CHIHARA-DUNKL-COEF}.
Upon imposing the condition~\eqref{SYM-COND}, one f\/inds the following equations for the terms in $\partial_{x}R$ and
$\partial_{x}$:
\begin{gather}
(x+\gamma)\omega(-x)+(-x+\gamma)\omega(x)=0,
\nonumber
\\
\omega'(x)=\left[\frac{\alpha}{x-\gamma}+\frac{\alpha+1}{x+\gamma}-\frac{2x\beta}{\gamma^2+1-x^2}\right]\omega(x).
\label{Eqs}
\end{gather}
It is easily seen that the common solution to~\eqref{Eqs} is given by
\begin{gather}
\label{CHIHARA-WEIGHT-2}
\omega(x)=\theta(x)(x+\gamma)\big(x^2-\gamma^2\big)^{\alpha}\big(1+\gamma^2-x^2\big)^{\beta},
\end{gather}
which corresponds to the weight function~\eqref{CHIHARA-WEIGHT-3} derived above.
It is directly checked that with~\eqref{CHIHARA-WEIGHT-2}, the equations for the terms in $\partial_{x}^{2}$, $R$ and
$\mathbb{I}$ arising from the symmetry condition~\eqref{SYM-COND} are identically satisf\/ied.
The orthogonality relation~\eqref{CHIHARA-ORTHO} can be recovered by the requirements that $\omega(x)>0$ on a~symmetric
interval.

\section[Chihara polynomials and big $q$ and $-1$ Jacobi polynomials]{Chihara polynomials and big $\boldsymbol{q}$ and $\boldsymbol{-1}$ Jacobi polynomials}\label{Section6}

In this section, the connexion between the Chihara polynomials and the big $q$-Jacobi and big~$-1$ Jacobi polynomials is
established.
In particular, it is shown that the Chihara polynomials are related to the former by a~$q\rightarrow -1$ limit and to
the latter by a~Christof\/fel transformation.

\subsection[Chihara polynomials and big $-1$ Jacobi polynomials]{Chihara polynomials and big $\boldsymbol{-1}$ Jacobi polynomials}

The big $-1$ Jacobi polynomials $J_{n}(x;a,b,c)$ were introduced in~\cite{Vinet-2012-05} as a~$q\rightarrow -1$ limit of the
big $q$-Jacobi polynomials.
In~\cite{Vinet-2011-02}, they were seen to be the polynomials that diagonalize the most general f\/irst order dif\/ferential
Dunkl operator preserving the space of polynomials of a~given degree (see~\eqref{-1-Jacobi-OP}).
The big $-1$ Jacobi polynomials can be def\/ined by their recurrence relation
\begin{gather}
\label{RECU-JACOBI}
xJ_{n}(x)=J_{n+1}(x)+(1-A_{n}-C_{n})J_{n}(x)+A_{n-1}C_{n}J_{n-1},
\end{gather}
subject to the initial conditions $J_{-1}(x)=0$, $J_{0}(x)=1$ and where the recurrence coef\/f\/icients read
\begin{gather}
\label{RECU-COEF-JACOBI}
A_{n}=
\begin{cases}
\dfrac{(1+c)(a+n+1)}{2n+a+b+2} & n \ \text{even},
\vspace{1mm}\\
\dfrac{(1-c)(n+a+b+1)}{2n+a+b+2} & n \ \text{odd},
\end{cases}
\qquad
C_{n}=
\begin{cases}
\dfrac{(1-c)n}{2n+a+b} & n \ \text{even},
\vspace{1mm}\\
\dfrac{(1+c)(n+b)}{2n+a+b} & n \ \text{odd},
\end{cases}
\end{gather}
for $0<c<1$.
Consider the monic polynomials $K_{n}(x)$ obtained from the big $-1$ Jacobi polynomials $J_{n}(x)$ by the Christof\/fel
transformation~\cite{Chihara-2011}
\begin{gather}
\label{Christoffel}
K_{n}(x)=\frac{1}{(x-1)}\left[J_{n+1}(x)-\frac{J_{n+1}(1)}{J_{n}(1)}J_{n}(x)\right]=(x-1)^{-1}[J_{n+1}(x)-A_{n}J_{n}(x)],
\end{gather}
where we have used the fact that
\begin{gather*}
J_{n+1}(1)/J_{n}(1)=A_{n},
\end{gather*}
which easily follows from~\eqref{RECU-JACOBI} by induction.
As is seen from~\eqref{Christoffel}, the polynomials $K_{n}(x)$ are kernel partners of the big $-1$ Jacobi polynomials
with kernel parameter $1$.
The inverse transformation, called the Geronimus transformation~\cite{Geronimus-1947}, is here given by
\begin{gather}
\label{Geronimus}
J_{n}(x)=K_{n}(x)-C_{n}K_{n-1}(x).
\end{gather}
Indeed, it is directly verif\/ied that upon substituting~\eqref{Christoffel} in~\eqref{Geronimus}, one recovers the
recurrence relation~\eqref{RECU-JACOBI} satisf\/ied by the big $-1$ Jacobi polynomials.
In the reverse, upon substituting~\eqref{Geronimus} in~\eqref{Christoffel}, one f\/inds that the kernel polynomials
$K_{n}(x)$ satisfy the recurrence relation
\begin{gather*}
xK_{n}(x)=K_{n+1}(x)+(1-A_{n}-C_{n+1})K_{n}(x)+A_{n}C_{n}K_{n-1}(x).
\end{gather*}
Using the expressions~\eqref{RECU-COEF-JACOBI} for the recurrence coef\/f\/icients, this recurrence relation can be cast in
the form
\begin{gather}
\label{recu}
x K_{n}(x)=K_{n+1}(x)+(-1)^{n+1}c K_{n}(x)+f_{n}K_{n-1}(x),
\end{gather}
where
\begin{gather}
\label{recu-coeff}
f_{n}=
\begin{cases}
\dfrac{(1-c^2)n(n+a+1)}{(2n+a+b)(2n+a+b+2)} & n \ \text{even},
\vspace{1mm}\\
\dfrac{(1-c^2)(n+b)(n+a+b+1)}{(2n+a+b)(2n+a+b+2)}& n \ \text{odd}.
\end{cases}
\end{gather}
It follows from the above recurrence relation that the kernel polynomials $K_{n}(x)$ of the big $-1$ Jacobi polynomials
correspond to the Chihara polynomials.
Indeed, taking $x\rightarrow x\sqrt{1-c^2}$ and def\/ining
\begin{gather*}
\alpha=b/2-1/2,
\qquad
\beta=a/2+1/2,
\qquad
\frac{c}{\sqrt{1-c^2}}=-\gamma,
\end{gather*}
it is directly checked that the recurrence relation~\eqref{recu} with coef\/f\/icients~\eqref{recu-coeff} corresponds to the
recurrence relation~\eqref{CHIHARA-RECU} satisf\/ied by the Chihara polynomials.
We have thus established that the Chihara polynomials are the kernel partners of the big $-1$ Jacobi polynomials with
kernel parameter $1$.
In view of the fact that the complementary Bannai--Ito polynomials are the kernel partners of the Bannai--Ito
polynomials, we have
\begin{gather*}
\begin{gathered}
\xymatrix
@C=.15\linewidth
@R=.1\linewidth
{
*+{
\begin{aligned}
&\text{Bannai--Ito}\\&B_{n}(x;\rho_1,\rho_2,r_1,r_2)
\end{aligned}}
\ar@<.5ex>[r]^-{\text{Christoffel}}
\ar@{->}[d]^-{h\rightarrow 0}_-{x\rightarrow x/h}
&
*+{\begin{aligned}
&\text{complementary BI}\\
&I_{n}(x;\rho_1,\rho_2,r_1,r_2)
\end{aligned}}
\ar@<.5ex>[l]^-{\text{Geronimus}}
\ar@{->}[d]^-{h\rightarrow 0}_-{x\rightarrow x/h}
\\
*+{
\begin{aligned}
&\text{big } {-}1\text{ Jacobi} \\
&J_{n}(x;a,b,c)
\end{aligned}
}
\ar@<.5ex>[r]^-{\text{Christoffel}}
&
*+{
\begin{aligned}
&\text{Chihara}\\
& C_{n}(x;\alpha,\beta,\gamma)
\end{aligned}
}
\ar@<.5ex>[l]^-{\text{Geronimus}}
}
\end{gathered}
\end{gather*}
The precise limit process from the Bannai--Ito polynomials to the big $q$-Jacobi polynomials can be found
in~\cite{Tsujimoto-2012-03}.

\subsection[Chihara polynomials and big $q$-Jacobi polynomials]{Chihara polynomials and big $\boldsymbol{q}$-Jacobi polynomials}

The Chihara polynomials also correspond to a~$q\rightarrow -1$ limit of the big $q$-Jacobi polynomials, dif\/ferent from the
one leading to the big $-1$ Jacobi polynomials.
Recall that the monic big $q$-Jacobi polynomials $P_{n}(x;\alpha,\beta,\gamma|q)$ obey the recurrence
relation~\cite{Koekoek-2010}
\begin{gather}
\label{RECU-Q-JACOBI}
x P_{n}(x)=P_{n+1}(x)+(1-\upsilon_{n}-\nu_{n})P_{n}(x)+\upsilon_{n-1}\nu_{n}P_{n-1}(x),
\end{gather}
with $P_{-1}(x)=0$, $P_{0}(x)=1$ and where
\begin{gather*}
\upsilon_{n}=\frac{(1-\alpha q^{n+1})(1-\alpha \beta q^{n+1})(1-\gamma q^{n+1})}{(1-\alpha \beta q^{2n+1})(1-\alpha
\beta q^{2n+2})},
\\
\nu_{n}=-\alpha \gamma q^{n+1}\frac{(1-q)^{n}(1-\alpha \beta \gamma^{-1}q^{n})(1-\beta q^{n})}{(1-\alpha \beta
q^{2n})(1-\alpha \beta q^{2n+1})}.
\end{gather*}
Upon writing
\begin{gather}
\label{Dompe}
q=-e^{\epsilon},
\qquad
\alpha=e^{2\epsilon \beta},
\qquad
\beta=-e^{\epsilon(2\alpha+1)},
\qquad
\gamma=-\gamma,
\end{gather}
and taking the limit as $\epsilon\rightarrow 0$, the recurrence relation~\eqref{RECU-Q-JACOBI} is directly seen to
converge, up to the redef\/inition of the variable $x\rightarrow x\sqrt{1-\gamma^2}$, to that of the Chihara
polynomials~\eqref{CHIHARA-RECU}.
In view of the well known limit of the Askey--Wilson polynomials to the big $q$-Jacobi polynomials, which can be found
in~\cite{Koekoek-2010}, we can thus write
\begin{gather*}
\begin{gathered}
\xymatrix
@C=.15\linewidth
@R=.08\linewidth
{
*+{
\begin{aligned}
&\text{Askey--Wilson}\\&p_{n}(x;a,b,c,d\,|\,q)
\end{aligned}}
\ar@{->}[r]^-{q\rightarrow -1}
\ar@{->}[d]^-{a\rightarrow 0}_-{x\rightarrow x/2a}
&
*+{\begin{aligned}
&\text{complementary BI}\\
&I_{n}(x;\rho_1,\rho_2,r_1,r_2)
\end{aligned}}
\ar@{->}[d]^-{h\rightarrow 0}_-{x\rightarrow x/h}
\\
*+{
\begin{aligned}
&\text{big } q\text{-Jacobi} \\
&P_{n}(x;\alpha,\beta,\gamma\,|\, q)
\end{aligned}
}
\ar@{->}[r]^-{q\rightarrow -1}
&
*+{
\begin{aligned}
&\text{Chihara}\\
& C_{n}(x;\alpha,\beta,\gamma)
\end{aligned}
}
}
\end{gathered}
\end{gather*}
It is worth mentioning here that the limit process~\eqref{Dompe} cannot be used to derive the bispectrality property of
the Chihara polynomials from the one of the big $q$-Jacobi polynomials.
Indeed, it can be checked that the $q$-dif\/ference operator diagonalized by the big $q$-Jacobi polynomials does not exist
in the limit~\eqref{Dompe}.
A similar situation occurs for the $q\rightarrow -1$ limit of the Askey--Wilson polynomials to the complementary
Bannai--Ito polynomials and is discussed in~\cite{Genest-2013-02-1}.

\section{Special cases and limits of Chihara polynomials}\label{Section7}

In this section, three special/limit cases of the Chihara polynomials $C_{n}(x;\alpha,\beta,\gamma)$ are considered.
One special case and one limit case correspond respectively to the generalized Gegenbauer and generalized Hermite
polynomials, which are well-known from the theory of symmetric orthogonal polynomials~\cite{Chihara-2011}.
The third limit case leads to a~new bispectral family of $-1$ orthogonal polynomials which depend on two parameters and
which can be seen as a~one-parameter extension of the generalized Hermite polynomials.

\subsection{Generalized Gegenbauer polynomials}

It is easy to see from the explicit expression~\eqref{CHIHARA-DEF} that if one takes $\gamma=0$, the Chihara polynomials
$C_{n}(x;\alpha,\beta,\gamma)$ become symmetric, i.e.\
$C_{n}(-x)=(-1)^{n}C_{n}(x)$.
Denoting by $G_{n}(x;\alpha,\beta)$ the polynomials obtained by specializing the Chihara polynomials to $\gamma=0$, one
directly has
\begin{gather*}
G_{2n}(x)=\frac{(-1)^{n}(\alpha+1)_{n}}{(n+\alpha+\beta+1)_{n}}\pFq{2}{1}{-n,n+\alpha+\beta+1}{\alpha+1}{x^2},
\\
G_{2n+1}(x)=\frac{(-1)^{n}(\alpha+2)_{n}}{(n+\alpha+\beta+2)_{n}}x\pFq{2}{1}{-n,n+\alpha+\beta+2}{\alpha+2}{x^2}.
\end{gather*}
The polynomials $G_{n}(x)$ are directly identif\/ied to the generalized Gegenbauer polynomials (see for
example~\cite{Belmehdi-2001-08,Dunkl-2001}).
In view of proposition~\eqref{Prop-1}, the polynomials $G_{n}(x)$ satisfy the recurrence relation
\begin{gather*}
xG_{n}(x)=G_{n+1}(x)+\sigma_{n}G_{n-1}(x),
\end{gather*}
with $G_{-1}(x)=0$, $G_{0}(x)=1$ and where $\sigma_{n}$ is given by~\eqref{CHIHARA-COEF}.
It follows from proposition~\eqref{Prop-2} that the polynomials $G_{n}(x)$ satisfy the family of eigenvalue equations
\begin{gather}
\label{ultra}
\mathcal{W}^{(\epsilon)}G_{n}(x)=\lambda^{(\epsilon)}_{n}G_{n}(x),
\end{gather}
where the eigenvalues are given by~\eqref{CHIHARA-EIGEN-2} and where the operator $\mathcal{W}^{(\epsilon)}$ has the
expression
\begin{gather*}
\mathcal{W}^{(\epsilon)}=S_{x}\partial_{x}^2+U_{x}\partial_{x}+V_{x}(\mathbb{I}-R),
\end{gather*}
with the coef\/f\/icients
\begin{gather*}
S_{x}=\frac{x^2-1}{4},
\qquad
U_x=\frac{x^2(\alpha+\beta+3/2)-\alpha-1/2}{2x},
\\
V_{x}=\frac{\alpha+1/2}{4x^2}-\frac{\alpha+\beta+3/2}{4}+\epsilon/2.
\end{gather*}
Upon taking
\begin{gather*}
\epsilon\rightarrow(\alpha+1)(\mu+1/2),
\qquad
\beta\rightarrow \alpha,
\qquad
\alpha\rightarrow \mu-1/2,
\end{gather*}
the eigenvalue equation~\eqref{ultra} can be rewritten as
\begin{gather}
\label{ultra-2}
\mathcal{Q}G_{n}(x)=\Upsilon_{n}G_{n}(x),
\end{gather}
where
\begin{gather}
\label{ultra-3}
\mathcal{Q}=(1-x^2)[D_{x}^{\mu}]^{2}-2(\alpha+1)xD_{x}^{\mu},
\end{gather}
where $D_{x}^{\mu}$ stands for the Dunkl derivative operator
\begin{gather}
\label{Dunkl-D}
D_{x}^{\mu}=\partial_{x}+\frac{\mu}{x}(\mathbb{I}-R),
\end{gather}
and where the eigenvalues $\Upsilon_{n}$ are of the form
\begin{gather}
\label{ultra-4}
\Upsilon_{2n}=-2n(2n+2\alpha+2\mu+1),
\qquad
\Upsilon_{2n+1}=-(2n+2\mu+1)(2n+2\alpha+2).
\end{gather}
The eigenvalue equation~\eqref{ultra-2} with the operator~\eqref{ultra-3} and eigenvalues~\eqref{ultra-4} corresponds to
the one obtained by Ben Cheikh and Gaied in their characterization of Dunkl-classical symmetric orthogonal
polynomials~\cite{BenCheikh-2007-04}.
The third proposition leads to the orthogonality relation
\begin{gather*}
\int_{-1}^{1}G_{n}(x)G_{m}(x)\omega(x)\mathrm{d}x=k_{n}\delta_{nm},
\end{gather*}
where the normalization factor $k_{n}$ is given by~\eqref{Norm} and where the weight function reads
\begin{gather*}
\omega(x)=|x|^{2\alpha+1}\big(1-x^2\big)^{\beta}.
\end{gather*}
We have thus established that the generalized Gegenbauer polynomials are $-1$ orthogonal polynomials which are
descendants of the complementary Bannai--Ito polynomials.

\subsection{A one-parameter extension of the generalized Hermite polynomials}

Another set of bispectral $-1$ orthogonal polynomials can be obtained upon letting
\begin{gather*}
x\rightarrow \beta^{-1/2}x,
\qquad
\alpha\rightarrow \mu-1/2,
\qquad
\gamma\rightarrow \beta^{-1/2}\gamma,
\end{gather*}
and taking the limit as $\beta\rightarrow \infty$.
This limit is analogous to the one taking the Jacobi polynomials into the Laguerre polynomials~\cite{Koekoek-2010}.
Let $Y_{n}(x;\mu,\gamma)$ denote the monic polynomials obtained from the Chihara polynomials in this limit.
The following properties of these polynomials can be derived by straightforward computations.
The polynomials $Y_{n}(x;\gamma)$ have the hypergeometric expression
\begin{gather*}
Y_{2n}(x)=(-1)^{n}(\mu+1/2)_{n}\pFq{1}{1}{-n}{\mu+1/2}{x^2-\gamma^2},
\\
Y_{2n+1}(x)=(-1)^{n}(\mu+3/2)_{n}(x-\gamma)\pFq{1}{1}{-n}{\mu+3/2}{x^2-\gamma^2}.
\end{gather*}
They satisfy the recurrence relation
\begin{gather*}
xY_{n}(x)=Y_{n+1}(x)+(-1)^{n}\gamma Y_{n}(x)+\vartheta_{n}Y_{n-1}(x),
\end{gather*}
with the coef\/f\/icients
\begin{gather*}
\vartheta_{2n}=n,
\qquad
\vartheta_{2n+1}=n+\mu+1/2.
\end{gather*}
The polynomials $Y_{n}(x;\gamma)$ obey the one-parameter family of eigenvalue equations
\begin{gather*}
\mathcal{Z}^{(\epsilon)}Y_{n}(x)=\lambda_{n}^{(\epsilon)}Y_{n}(x),
\end{gather*}
where the spectrum has the form
\begin{gather*}
\lambda_{2n}^{(\epsilon)}=n,
\qquad
\lambda_{2n+1}^{(\epsilon)}=n+\epsilon.
\end{gather*}
The explicit expression for the second-order dif\/ferential Dunkl operator $\mathcal{Z}^{(\epsilon)}$ is
\begin{gather*}
\mathcal{Z}^{(\epsilon)}=S_{x}\partial_{x}^2-T_{x}\partial_{x}R+U_{x}\partial_{x}+V_{x}(\mathbb{I}-R),
\end{gather*}
with
\begin{gather*}
S_{x}=\frac{\gamma^2-x^2}{4x^2},
\qquad
T_{x}=\frac{\gamma(x-\gamma)}{4x^3},
\\
U_x=\frac{x}{2}+\frac{\gamma}{4x^2}-\frac{\gamma^2}{2x^3}-\frac{\mu+\gamma^2}{2x},
\qquad
V_{x}=\frac{3\gamma^2}{8x^{4}}-\frac{\gamma}{4x^3}+\frac{\mu+\gamma^2}{4x^2}+\epsilon\frac{x-\gamma}{2x}-\frac{1}{4}.
\end{gather*}
The algebra encoding this bispectrality of the polynomials $Y_{n}(x)$ is obtained by taking
\begin{gather*}
K_1=\mathcal{Z}^{(\epsilon)},
\qquad
K_2=x,
\qquad
P=R+\frac{\gamma}{x}(\mathbb{I}-R),
\end{gather*}
and def\/ining $K_3=[K_1,K_2]$.
One then has the commutation relations
\begin{gather*}
[K_1,P]=0,
\qquad
\{K_2,P\}=2\gamma,
\qquad
\{K_3,P\}=0,\\
[K_2,K_3]=(2\epsilon-1)K_2^2P-2\gamma K_3P+(\gamma^2-2\gamma\epsilon+\mu)P+1/2,
\\
[K_3,K_1]=(1-2\epsilon)K_3P+\epsilon(\epsilon-1)K_2+\gamma\epsilon(\epsilon-1)P.
\end{gather*}
The orthogonality relation reads
\begin{gather*}
\int_{\mathcal{S}}Y_{n}(x)Y_{m}(x)w(x)\mathrm{d}x=l_{n}\delta_{nm}
\end{gather*}
with $\mathcal{S}=(-\infty,-|\gamma|]\cup[|\gamma|,\infty)$ and with the weight function
\begin{gather*}
w(x)=\theta(x)(x+\gamma)\big(x^2-\gamma^2\big)^{\mu-1/2}e^{-x^2}.
\end{gather*}
The normalization factors
\begin{gather*}
l_{2n}=n!e^{-\gamma^2}\Gamma(n+\mu+1/2),
\qquad
l_{2n+1}=n!e^{-\gamma^2}\Gamma(n+\mu+3/2)
\end{gather*}
are obtained using the observation that the polynomials $Y_{n}(x;\gamma)$ can be expressed in terms of the standard
Laguerre polynomials~\cite{Koekoek-2010}.

\subsection{Generalized Hermite polynomials}

The polynomials $Y_{n}(x;\mu,\gamma)$ can be can be considered as a~one-parameter extension of the ge\-ne\-ralized Hermite
polynomials.
Indeed, upon denoting by $H_{n}^{\mu}(x)$ the polynomials obtained by taking $\gamma=0$ in $Y_{n}(x;\mu,\gamma)$, one
f\/inds that
\begin{gather*}
H_{2n}^{\mu}(x)=(-1)^{n}(\mu+1/2)_{n}\pFq{1}{1}{-n}{\mu+1/2}{x^2},
\\
H_{2n+1}^{\mu}(x)=(-1)^{n}(\mu+3/2)_{n}x\pFq{1}{1}{-n}{\mu+3/2}{x^2},
\end{gather*}
which corresponds to the generalized Hermite polynomials~\cite{Chihara-2011}.
It is thus seen that the ge\-ne\-ralized Hermite polynomials are also $-1$ orthogonal polynomials that can be obtained from
the complementary Bannai--Ito polynomials.
For this special case, the eigenvalue equations can be written (taking $\epsilon\rightarrow \epsilon/2$) as
\begin{gather*}
\Omega^{(\epsilon)}H_{n}^{\mu}(x)=\lambda_{n}^{(\epsilon)}H_{n}(x),
\qquad \mbox{with}
\qquad
\lambda_{2n}^{(\epsilon)}=2n,
\qquad
\lambda_{2n+1}^{(\epsilon)}=2n+\epsilon
\end{gather*}
and where the operator $\Omega^{(\epsilon)}$ reads
\begin{gather*}
\Omega^{(\epsilon)}=-\frac{1}{2}\partial_{x}^2+\left(x-\frac{\mu}{x}\right)\partial_{x}+\left(\frac{\mu}{2x^2}+\frac{\epsilon-1}{2}\right)(\mathbb{I}-R).
\end{gather*}
The orthogonality relation then reduces to
\begin{gather*}
\int_{-\infty}^{\infty}H_{n}^{\mu}(x)H_{m}^{\mu}(x)|x|^{\mu}e^{-x^2}\mathrm{d}x=l_{n}\delta_{nm}.
\end{gather*}
Upon taking $\widetilde{\Omega}^{(\epsilon)}=e^{-x^2/2}\Omega^{(\epsilon)}e^{x^2/2}$, the eigenvalue equations can be
written as
\begin{gather*}
\widetilde{\Omega}^{(\epsilon)}\psi_{n}(x)=\lambda_{n}^{(\epsilon)}\psi_{n}(x),
\end{gather*}
with $\psi_{n}=e^{-x^2/2}H_{n}^{\mu}(x)$ and with the eigenvalues
\begin{gather*}
\lambda_{2n}^{(\epsilon)}=2n+\mu+1/2,
\qquad
\lambda_{2n+1}^{(\epsilon)}=2n+\mu+3/2+\epsilon.
\end{gather*}
The operator $\widetilde{\Omega}^{(\epsilon)}$ can be cast in the form
\begin{gather*}
\widetilde{\Omega}^{(\epsilon)}=-\frac{1}{2}(D_{x}^{\mu})^2+\frac{1}{2}x^2+\frac{\epsilon}{2}(\mathbb{I}-R),
\end{gather*}
where $D_{x}^{\mu}$ is the Dunkl derivative~\eqref{Dunkl-D}.
The operator $\widetilde{\Omega}^{(\epsilon)}$ corresponds to the Hamiltonian of the one-dimensional Dunkl
oscillator~\cite{Rosenblum-1994}.
Two-dimensional versions of this oscillator models have been considered recently~\cite{Genest-2013-04,Genest-2013-09,Genest-2013-07}.

\section{Conclusion}\label{Section8}

In this paper, we have characterized a~novel family of $-1$ orthogonal polynomials in a~con\-tinuous variable which are
obtained from the complementary Bannai--Ito polynomials by a~limit process.
These polynomials have been called the Chihara polynomials and it was shown that they diagonalize a~second-order
dif\/ferential Dunkl operator with a~quadratic spectrum.
The orthogonality weight, the recurrence relation and the explicit expression in terms of Gauss hypergeometric function
have also been obtained.
Moreover, special cases and descendants of these Chihara polynomials have been examined.
From these considerations, it was observed that the well-known generalized Gegenbauer/Hermite polynomials are in fact~$-1$ polynomials.
In addition, a~new class of bispectral~$-1$ orthogonal polynomials which can be interpreted as a~one-parameter extension
of the generalized Hermite polynomials has been def\/ined.

With the results presented here, the polynomials in the higher portion of the emerging tableau of $-1$ orthogonal
polynomials are now identif\/ied and characterized.
At the top level of the tableau sit the Bannai--Ito polynomials and their kernel partners, the complementary Bannai--Ito
polynomials.
Both sets depend on four parameters.
At the next level of this~$-1$ tableau, with three parameters, one has the big~$-1$ Jacobi polynomials, which are
descendants of the BI polynomials, as well as the dual~$-1$ Hahn polynomials (see~\cite{Tsujimoto-2013-03-01}) and the
Chihara polynomials which are descendants of the CBI polynomials.
The complete tableau of~$-1$ polynomials with arrows relating them shall be presented in an upcoming review.

\subsection*{Acknowledgements}

V.X.G.\ holds a~fellowship from the Natural Sciences and Engineering Research Council of Canada (NSERC).
The research of L.V.\ is supported in part by NSERC.

\pdfbookmark[1]{References}{ref}
\LastPageEnding

\end{document}